\magnification=\magstep1
\input amssym.def 
\input amssym
\hsize=5.25in
\vsize=7in
\hfill{Ref. SISSA 42/2000/FM}
\vskip2truecm

\centerline{
{\bf $A(SL_q(2))$ at roots of unity is a free module over $A(SL(2))$}}

\vskip2truecm
\centerline{Ludwik D\c abrowski, Cesare Reina and Alessandro Zampa}
\vskip.2truecm
\centerline{SISSA, via Beirut 2-4, 34014 Trieste (Italy)}

\centerline{dabrow@sissa.it, reina@sissa.it, zampa@fm.sissa.it}

\vskip2truecm
\hsize=4.25in
{\leftskip=.5in
\noindent
{{\bf Abstract.} It is shown that when $q$ is a primitive root of unity
of order not equal to $2$ mod $4$,
$A(SL_q(2))$ is a {\it free} module of finite rank over the coordinate ring of
the classical group $SL(2)$. An explicit set of generators is
provided.
}
\bigskip\par}
\vskip.3truecm
\centerline{~~MSC: ~17B37, 16W35, 81R50}
\centerline{~Key words: quantum group}

\vskip3.2truecm
\hsize=5.25in
\hoffset=.0in
Quantum groups at roots of unity have particularly rich and interesting structure.
In this letter we adopt the `quantum functions' point of view, 
dual to the quantum universal enveloping algebras (cf. [7]).
It is known that the quantum Hopf algebra $A(SL_q(2))$ at primitive roots of unity 
$q^l = 1, ~l$ odd, or $q^{2l} = 1, ~l$ even,  
can be viewed as a module over the algebra $A(SL(2))$,
eg. by using the Frobenius mappings [9, 11].
This module is projective and finitely generated [5] and thus
it corresponds to a locally free coherent sheaf $\cal F$ over 
the affine group scheme $SL(2)=:SL(2,{\Bbb C})$. 
It has been also established that $SL_q(2)$ at primitive odd roots of unity 
forms a principal fibre bundle (faithfully flat Hopf-Galois extension) 
over $SL(2)$ [1, 8]. (See also [3], [4], [6] for the case of cubic root of unity, 
that was conjectured in [2] to be important in understanding the (quantum) 
symmetry of the Standard Model of fundamental interactions). 

It may have been guessed for a long time that 
$A(SL_q(2))$ at roots of unity is actually free over $A(SL(2))$ but,  
as far as we know, no proof has been given. 
Here we show that it is actually so. Our reasoning 
underlines the difficulties to tackle the more general problem of
other ``simple'' quantum groups at roots of unity.  

Recall that for ${\Bbb C}\ni q \neq 0$, $A(SL_q(2))$ is the free algebra 
${\Bbb C}\!\! <\! a,b,c,d\! >$
modulo the ideal generated by the commutation relations
$$ab=qba\quad ac=qca,\quad bd=qdb,\quad bc=cb,\quad cd=qdc,\quad
 ad-da=(q-q^{-1})bc$$
and by the $q$-determinant relation $ad-qbc=1.$
When $q$ is a primitive $l$-th root of unity, for odd $l$, 
the Hopf subalgebra generated by 
$$\alpha = a^l,\; \beta=b^l,\; \gamma= c^l,\; \delta=d^l$$
is isomorphic to the Hopf algebra 
$$A(SL(2))={\Bbb C}[\alpha ,\beta ,\gamma , \delta ]\; /
<\alpha\delta - \beta\gamma -1>$$
with the restricted coproduct, counit and coinverse.
Moreover, $A(SL_q(2))$ is a finitely generated module over $A(SL(2))$.
The same happens when $l$ is even and $q$ is a primitive $2l$-th root of unity, 
but this time $A(SL_q(2))$ has different left and right module structures, 
due to the appearance of some $\pm$ signs in the commutation relations. 

Notice that, using the commutation relations to order
monomials, $A(SL_q(2))$ can be generated over $A(SL(2))$ by the 
$l^4$ products 
$$a^{r_a}b^{r_b}c^{r_c}d^{r_d}, \quad (0\leq r_\cdot\leq l-1).$$
The determinant relation $ad=1+qbc$ reduces the generators above
to a smaller set; all the monomials
with $r_ar_d>0$ can be expressed as linear combinations over 
${\Bbb C}$ of the following set
  $$a^{r_a}b^{r_b}c^{r_c},\quad b^{r_b}c^{r_c}d^{r_d},$$
of $2l^3$ elements. Indeed
 $$P_k(bc)=:a^kd^k=\prod_{j=1}^k\; (1+q^{2j-1}bc)=
  \sum_{j=0}^kp_{k,j}(q)b^jc^j,$$
where 
$$p_{k,j}(q)=q^{j^2} \prod_{r=j+1}^{k} (q^{2r} - 1) 
~~\Bigl(~\prod_{s=1}^{k-j} (q^{2s} - 1)\Bigr)^{-1} .$$
 
Since $A(SL_q(2))$ is a finitely generated projective module 
over $A(SL(2))$, it corresponds to a locally free coherent 
sheaf $\cal F$ over $Spec(A(SL(2)))=SL(2)$.
We now think of $SL(2)$ as a complex
manifold $SL(2)_{an}$ instead as an affine scheme.

\vskip3truemm\noindent
{\bf Proposition 1.} {\it ${\cal F}\rightarrow SL(2)_{an}$ is 
locally free of rank $l^3$ and 
the corresponding vector bundle $F$ is trivial.}

\vskip3truemm\noindent
{\bf Proof.} We cover $SL(2)_{an}$ with two open sets $U_\alpha,\; U_\beta$
where $\alpha$ or $\beta$ is invertible. Since $a^{-1}=\alpha^{-1}a^{l-1}$, 
also $a^{-1}$ exists on $U_\alpha$ and similarly $b^{-1}$ on $U_\beta$.
Thus on $U_\alpha$ we can set $d=a^{-1}(1+qbc)$ and use
the $l^3$ elements $a^{r_a}b^{r_b}c^{r_c}$ as generators, 
and on $U_\beta$ we set $c=q^{-1}b^{-1}(ad-1)$
and use $a^{r_a}b^{r_b}d^{r_d}$ as generators.

\noindent
This proves that $\cal F$ is locally free of rank $l^3$
and hence it is a sheaf of sections of a vector bundle $F$ over  
$SL(2)_{an}$.
But $SL(2)_{an}$ retracts to its maximal compact subgroup $SU(2)\simeq S^3$
and since for every Lie group $G$, $\pi_2 (G) =0$, $F$ is topologically trivial.
Now $SL(2)_{an}$ is a Stein domain, the fibre bundles 
over it have no moduli and are simply classified topologically.
It follows that $F$ is also holomorphically trivial.
\hfill{$\bullet$}

\vskip3truemm
\noindent
The triviality of $F$ (and of $\cal F$) implies that we should be
able to find $l^3$ global nowhere vanishing independent sections.
It is clear that the local trivializations given above are by no 
means global. 
To find out global sections of $F$ we work algebraically. \
The following theorem provides an explicit
set of $l^3$ generators of $A(SL_q(2))$ over $A(SL(2))$.

\vskip3truemm\noindent
{\bf Proposition 2.} {\it Let $q$ be a primitive 
root of unity of order $l$ with $l$ odd,
or of order $2l$ with $l$ even. Then
$A(SL_q(2))$ is a free (left and/or right) module of rank 
$l^3$ over the coordinate ring of the classical group $A(SL(2))$.
The $l^3$ generators can be chosen as
$$a^mb^nc^{s'},\qquad b^nc^{s''}d^r,$$
with the integers $m, n, r, s', s''$ in the range
$1\le m\le l\! -\! 1;\; 0\le n,r\le l\! -\! 1;\; m\le s'\le l\! -\! 1;\; 
0\le s''\le l\! -\! r\! -\! 1$.}
         
\vskip3truemm\noindent
{\bf Proof.} For concreteness, we choose the left module structure,
the right module case is analogous.
 Using the above expression for $P_k(bc)$ and 
noticing that $p_{k,0}(q)=1$, one has
$$\eqalign{
  b^rc^sd^{l-k}=&~
  q^{(l-k)(r+s)} ~\delta\; a^k b^r c^s -
  \sum _{j=1}^k p_{k,j}(q)~ b^{r+j}c^{s+j}d^{l-k},\cr 
  a^{l-k}b^rc^s=&~
  q^{-k(r+s)} ~\alpha\; b^r c^s d^k -
  \sum _{j=1}^k p_{k,j}(q)~ a^{l-k}b^{r+j}c^{s+j}.\cr}$$
The idea is to use recursively these relations to eliminate the 
monomials in the l.h.s.'s in terms of those in the r.h.s.'s. Of
course there is a cyclicity problem, because the recursive
substitutions will sooner or later bring in the r.h.s the monomial
in the l.h.s. multiplied by a non constant element in $A(SL(2))$.
We then have to put bounds on the ranges of the exponents.
Let us look at the first equation. As $d$ occurs with the same 
exponent in both sides, there are no restrictions on $k$. 
Fix $s=s_1$ (possibly depending on $k$), 
let $S_1(k)=\{(s_1+j)\bmod l,\; 0\le j\le k\}$,
$S_1^c(k)={\Bbb Z}_l- S_1(k)$.
We can use the monomials
$b^rc^sd^{l-k-1}$ with $s\in S_1(k)$, $k,r\in {\Bbb Z}_l$ and 
$a^kb^rc^s$ with $k,r,s\in {\Bbb Z}_l$ to generate the monomials 
in the l.h.s. with $s\in S_1^c(k)$. The same can be done for the
second equation, by choosing $s_2$. We get a system of generators
of the form $b^rc^sd^{l-k-1}$, with $r,k\in {\Bbb Z}_l$, $s\in S_1(k)$,
and $a^{k}b^rc^s$, with $k\in {\Bbb Z}_l-\{ 0\}$, $r\in {\Bbb Z}_l$ and
$s\in S_2(l-k-1)$. These are $l^3$ elements and
we have to choose carefully $s_1, s_2$ to be sure that there are no
relations among them. This is the case for $s_1=0$ and $s_2=k$. 
By a suitable relabelling of the indices, we complete the proof of the theorem.
\hfill{$\bullet$}

\vskip3truemm
\noindent {\bf Remark 1.} The results above hold for 
$q$ being a primitive root of unity of order not equal to $2$ mod $4$.
In the remaining cases, $q^{2l} = 1$, $l$ odd, 
the $2l$-th powers of $a,\; b,\; c,\; d$
commute but neither the determinant relation nor
the coproduct close on the algebra generated by them. 
To remedy this mismatch one can enlarge the algebra by taking 
$\alpha ,\; \beta ,\; \gamma ,\; \delta$
(the $l$-th powers of $a,\; b,\; c,\; d$) as generators, 
obtaining however a non-commutative algebra.

\vskip3truemm
\noindent {\bf Remark 2.} In attempting to generalize the above results
to $SL_q(n)$ we immediately face the lack of the 
topological argument used in the proof of Proposition 1, 
which in fact has been the starting point of this letter. 
So in principle one has to use algebra alone.
Local triviality easily follows by considering the open covering
of $SL(n)$ given by the subsets $U_k$ where the algebraic complements
$\hat T_{k1}$ of the elements of the first column are invertible
(see [10] for notations). On $U_k$ we have   
$$T_{k1}=(-q)^{1-k}
 [1-\sum_{j\ne k}(-q)^{j-1} T_{j1}\hat T_{j1}]\hat T_{k1}^{-1}$$
and a basis of the localized module is given by the $l^{n^2-1}$
monomials $\prod_{(r,s)\ne (k,1)} T_{rs}^{n_{rs}}$, with 
$n_{rs}\in {\Bbb Z}_l$. This simplifies in the special case of $SL_q(n)$
the proof in [5] valid for a general quantum group.
The determinant relation provides equations similar to those of Proposition 2,
but the elimination problem is much subtler. 

\vskip.5truecm\noindent
{\bf Bibliography}
\vskip.3truecm\noindent
\item{1.}{N. Andruskiewitsch, 
{\it Notes on Extensions of Hopf algebras}
Can. J. Math. {\bf 48} 3-42 (1996)}  

\item{2.}{A. Connes,
{\it Gravity coupled with matter and the foundation of non commutative geometry},
Commun.~Math.~Phys. {\bf 182}  155--176 (1996)}

\item{3.}{R. Coquereaux, A. O. Garcia, R. Trinchero,
{\it Hopf stars, twisted Hopf stars and scalar products on quantum spaces},
math-ph/9904037} 

\item{4.}{L. D\c abrowski, P. M. Hajac, P. Siniscalco,~ 
{\it Explicit Hopf-Galois Description of $SL\sb{e\sp{2\!\pi\! i\! /\!
3}}(2)$-Induced Frobenius Homomorphisms}, in: Enlarged Proceedings of
the ISI GUCCIA Workshop on quantum groups, non commutative
geometry and fundamental physical interactions, 
    1997 Palermo, Italy; D.Kastler, M.Rosso, T.Schucker eds. 
    Nova Science Pub. Inc., Commack, New-York 279--299, 1999}

\item{5.}{C. De Concini, V. Lyubashenko, 
{\it Quantum function algebra at roots of 1},
Adv. Math. {\bf 108}  205-262 (1994)}

\item{6.}{D. Kastler, {\it $U_q(sl_2)$ for $|q|\! =\! 1$
as the complexification of a real Hopf algebra},
Nachrichten der Academie der Wissenschaften in Gottingen,
Vanderhoeck and Ruprecht, 1999}

\item{7.}{G. Lusztig,
{\it Quantum Groups at Roots of Unity}, Geom. Dedicata {\bf 35} 89--114 (1991)}

\item{8.}{S. Montgomery, H.J. Schneider,
{\it Prime ideals in Hopf Galois extensions}, in prep.}

\item{9.}{B. Parshall, J. Wang,
{\it Quantum linear groups.} 
American Mathematical Society Memoirs no. 439,
Providence, R.I., American Mathematical Society (AMS) (1991)}

\item{10.}{C. Reina, A. Zampa, 
{\it Quantum homogeneous spaces at roots of unity}, 
in: Quantization, Coherent States and Poisson Structures,
A. Strasburger et al. eds, Polish Sci. Pub. PWN - Warszawa 1998}

\item{11.}{M. Takeuchi,
{\it Some topics on $GL_q(n)$},
J. Alg. {\bf 147} 379-410 (1992)}

\bye